\documentclass{amsart}

\usepackage{amssymb,amsmath,amsthm,latexsym}
\usepackage{epsfig}
\usepackage{psfrag}
\newtheorem{Theorem}{\bf Theorem}
\newtheorem{Lemma}[Theorem]{\bf Lemma}
\newtheorem{Proposition}[Theorem]{\bf Proposition}
\newtheorem{Corollary}[Theorem]{\bf Corollary}

\def\qed{\hfill$\Box$}
\newcommand{\be}{\begin{equation}}
\newcommand{\ee}{\end{equation}}

\def\qed{\hfill$\Box$}

\def\C{\mathbb C}

\def\k{{\bf k}}

\def\CE{{\mathcal E}}

\def\hpic #1 #2 {\mbox{$\begin{array}[c]{l} 
\epsfig{file=#1,height=#2}\end{array}$}}
\def\wpic #1 #2 {\mbox{$\begin{array}[c]{l} 
\epsfig{file=#1,width=#2}\end{array}$}}

\begin{document}

\title[The planar algebra of a semisimple and cosemisimple Hopf algebra]{The 
planar algebra of a semisimple and cosemisimple Hopf algebra}

\author{Vijay Kodiyalam}
\address{The Institute of Mathematical Sciences, Chennai, India}
\email{vijay@imsc.res.in,sunder@imsc.res.in}
\author{V. S. Sunder}

\subjclass{Primary 16W30; Secondary 46L37}
%\subjclass{Primary 17B20, 13A50; Secondary 57M27}
%\date{January 1, 1994 and, in revised form, June 22, 1994.}

%\keywords{Differential geometry, algebraic geometry}

\begin{abstract}
To a semisimple and cosemisimple Hopf algebra over an algebraically closed 
field,
we associate a planar algebra defined by generators and relations and
show that it is a connected, irreducible, spherical, non-degenerate planar 
algebra
with non-zero modulus and 
of depth two.
This association is shown to yield a bijection between (the isomorphism
classes, on both sides, of) such objects. 
\end{abstract}

\maketitle

Throughout this paper, the symbol $\k$ will always denote an
algebraically closed field and $H = (H,\mu,\eta,\Delta,\epsilon,S)$ will 
always denote a semisimple and
cosemisimple (necessarily finite-dimensional) Hopf algebra over $\k$. 
We associate to $H$, a planar algebra over the field $\k$
which is an analogue of the construction in \cite{KdyLndSnd} of
the `subfactor planar
algebra' associated to a (finite-dimensional) Kac algebra.

We then study various properties of this planar algebra including
computation of its partition function and duality with the planar
algebra of $H^*$.
Conversely, we show that every connected, irreducible, spherical,
non-degenerate planar algebra with non-zero modulus and of depth two
arises in this manner, thus obtaining a generalisation of the
Ocneanu-Szymanski theorem (cf. \cite{Szy}).

\section{Semisimple and cosemisimple Hopf algebras}
We begin by recalling well-known facts about such an $H$, the proofs of which
may be found in \cite{TngGlk} and in \cite{LrsRdf}.
The semisimplicity and cosemisimplicity assumptions imply that
both $H$ and $H^*$ are multi-matrix algebras and  the 
dimensions, say $n$, of $H$, as well as those of its irreducible 
representations, are `non-zero in $\k$'.
It follows that the traces
in the regular representations of $H$ and $H^*$, which we shall denote by 
$\phi$
and $h$ respectively, are non-degenerate traces.
Further, these are two-sided integrals for $H^*$ and $H$ respectively : i.e., 
they satisfy $\phi\psi = 
\psi(1)\phi = \psi\phi$ and
$hx = \epsilon(x)h = xh$, for all $\psi \in H^*$ and $x \in H$.
Also, $\epsilon(h) = \phi(h) = \phi(1) = n \in \k$.
Finally, the antipodes of $H$ and $H^*$ are involutive.

We shall use the standard notations $\mu_k$ and $\Delta_k$ for the
$k$-fold iterated product and coproduct respectively.
In particular, $\mu_0 = \eta$, $\mu_1 = id_H$, $\mu_2 = \mu$ and
$\Delta_0 = \epsilon$, $\Delta_1 = id_H$ and $\Delta_2 = \Delta$. 
We will find it convenient to use our version of the so-called Sweedler
notation for comultiplication - according to which we write, for example,
$\Delta_n(x) = x_1 \otimes \cdots \otimes x_n$
rather than the more familiar
$\Delta_n(x) = \sum_{(x)}x_{(1)} \otimes \cdots \otimes x_{(n)}$
in the interest of notational  convenience.

\section{Planar algebras}
We will also need the formalism of Jones' planar algebras.
Although Jones primarily used planar algebras which are $C^*$-planar
algebras, and {\em \`{a} fortiori} defined over $\C$, we need their
analogues over more general fields here, so we give a `crash course'
accordingly; but by and large, we assume familiarity with planar
tangles, planar networks and planar 
algebras. The basic reference is \cite{Jns}. A somewhat more
leisurely treatment of the basic notions may also be found in \cite{KdySnd}.
(We will follow the latter where, for instance, the $*$'s are attached to
`distinguished points' on boxes rather than to regions and for notation
such as $1_k$ for the identity element of $P_k$ for a planar algebra $P$.)

We shall continue to use the term {\em planar tangle}, as well its
{\em colour}, in exactly the
same sense as used in \cite{KdySnd}. By a planar algebra $P$ (over
$\k$), we shall mean a collection $\{P_k: k \in Col\}$ of $\k$-vector
spaces, indexed by
the set $Col = \{0_+, 0_-, 1, 2, \cdots \}$ of `colours', which comes
equipped with the following structure: to each $k_0$-tangle $T$ (with
internal boxes $B_1, \cdots, B_g$ of colours $k_1, \cdots , k_g$) is
associated a $\k$-linear mapping $Z^P_T:\otimes_{i=1}^g P_{k_i}
\rightarrow P_{k_0}$ which satisfies several natural conditions
(listed as equations (2.2), (2.3), (2.4) and (2.6) in \cite{KdySnd}).

Given a `label set' $L = \coprod_{k \in Col} L_k$, an {\em L-labelled
  tangle} is a tangle $T$ equipped with a labelling of every internal
  box of colour $k$ by an element from $L_k$. The {\em
  universal planar algebra on $L$}, denoted by $P(L) =
  \{P(L)_k: k \in ~Col\}$ is defined by requiring that  $P(L)_k$ is the
  $\k$-vector space with basis consisting of the set of all
  $L$-labelled $k$-tangles, with  the action of a planar
  tangle on a tensor product of basis vectors given by $L$-labelled
  tangles being  the obvious $L$-labelled tangle obtained by substitution.

Recall that a planar ideal of a planar algebra $P$ is a set $I = \{I_k:k
\in ~Col\}$ with the property that (i) each $I_k$ is a subspace of
$P_k$, and (ii) for any $k_0$-tangle $T$ as before, $Z_T(\otimes_{i=1}^g x_i) 
\in
I_{k_0}$ whenever $x_i \in I_{k_i}$ for at least one
$i$. Given  a planar ideal $I$ in a planar algebra $P$, there is a
natural planar algebra structure on the `quotient' $P/I = \{P_k/I_k :
k \in  ~Col\}$. 

Given any subset $R=\{R_k: k \in ~Col\}$ of $P$ (meaning $R_k \subset
P_k$), there is a smallest planar ideal $I(R) = \{I(R)_k:k \in ~Col\}$
such that $R_k \subset I(R)_k$ for all $k \in Col$; and finally, given a label
set $L$ as above, and any 
subset $R$ of the universal planar algebra $P(L)$, the quotient
$P(L)/I(R)$ is said to be the planar algebra $P(L,R)$ {\em
presented with generators $L$ and relations $R$}.

\section{Definition of $P(H,\delta)$}
Motivated by the results of \cite{KdyLndSnd} - where the case of the
so-called finite-dimensional Kac-algebras (over $\C$) is treated - we
wish to define the planar 
algebra associated to  a semisimple and cosemisimple Hopf
algebra via generators and relations.

\begin{figure}[!h]
\begin{center}
\psfrag{zab}{\huge $\zeta a + b$}
\psfrag{eq}{\huge $=$}
\psfrag{a}{\huge $a$}
\psfrag{b}{\huge $b$}
\psfrag{z}{\huge $\zeta$}
\psfrag{+}{\huge $+$}
\psfrag{del}{\huge $\delta$}
\resizebox{12.0cm}{!}{\includegraphics{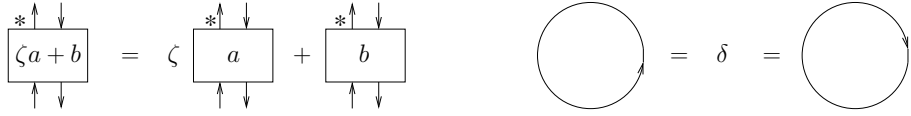}}
\end{center}
\caption{The L(inearity) and M(odulus) relations}
\label{fig:LnrMdl}
\end{figure}

\begin{figure}[!h]
\begin{center}
\psfrag{zab}{\huge $\zeta a + b$}
\psfrag{eq}{\huge $=$}
\psfrag{1h}{\huge $1_H$}
\psfrag{h}{\huge $h$}
\psfrag{z}{\huge $\zeta$}
\psfrag{+}{\huge $+$}
\psfrag{del}{\huge $\delta$}
\resizebox{12.0cm}{!}{\includegraphics{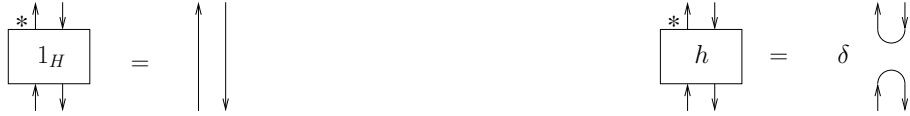}}
\end{center}
\caption{The U(nit) and I(ntegral) relations}
\label{fig:NitNtg}
\end{figure}

\begin{figure}[!h]
\begin{center}
\psfrag{epa}{\huge $\epsilon(a)$}
\psfrag{eq}{\huge $=$}
\psfrag{delinphia}{\huge $\delta^{-1} \phi(a)$}
\psfrag{h}{\huge $h$}
\psfrag{a}{\huge $a$}
\psfrag{+}{\huge $+$}
\psfrag{del}{\huge $\delta$}
\resizebox{12.0cm}{!}{\includegraphics{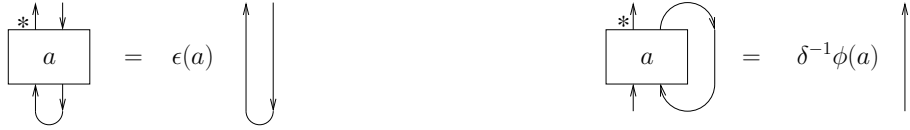}}
\end{center}
\caption{The C(ounit) and T(race) relations}
\label{fig:CntTrc}
\end{figure}

\begin{figure}[!h]
\begin{center}
\psfrag{epa}{\huge $\epsilon(a)$}
\psfrag{eq}{\huge $=$}
\psfrag{delinphia}{\huge $\delta^{-1} \phi(a)$}
\psfrag{a1}{\huge $a_1$}
\psfrag{a2b}{\huge $a_2\,b$}
\psfrag{b}{\huge $b$}
\psfrag{a}{\huge $a$}
\psfrag{sa}{\huge $Sa$}
\psfrag{del}{\huge $\delta$}
\resizebox{12.0cm}{!}{\includegraphics{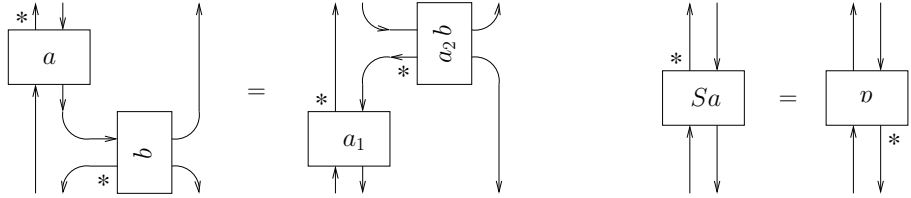}}
\end{center}
\caption{The E(xchange) and A(ntipode) relations}
\label{fig:XchNtp}
\end{figure}

However, our definition will depend on a choice we have to make of a
square root in $\k$ of $n$. To be precise, we shall let $\delta$ be a solution
to $\delta^2 = n$ in $\k$, and then define $P(H,\delta)$
to be the planar algebra $P(L,R)$, with
\[ L_k  = \left\{ \begin{array}{ll} H & {\mbox {if}} ~ k=2\\\emptyset &
  {\mbox{otherwise}} \end{array} \right. \]
and ~$R$ being given by the set of relations in Figures \ref{fig:LnrMdl} - 
\ref{fig:XchNtp} (where (i) we write 
the relations as identities - so the statement 
$a=b$ is interpreted as $a-b \in R$; (ii) $\zeta \in \k$  and  $a,b \in H$; and
(iii) the external boxes of all tangles appearing in the relations are left 
undrawn and it is assumed that all external $*$'s are at the top left corners).

We note that relation (M) actually represents two relations, one in $P_{0_+}$
and the other in $P_{0_-}$ and that the  $\delta$ in relation (M) means $\delta 1_{0_+}$ in
one of the relations and $\delta 1_{0_-}$ in the other.

In the rest of this paper, we shall simply write $P$ for the planar
algebra $P(H,\delta)$.

\section{Properties of $P$}
We wish to study the properties of $P$.
Recall that a planar algebra $P$ is said to be connected if $dim~P_{0_\pm} = 
1$ and irreducible if $dim ~P_1 = 1$.
A connected planar algebra is said to have modulus $\delta$ if the relation 
(M) holds in $P$.

We will need some facts concerning `exchange relation planar algebras'
as defined in \cite{Lnd}. Although only planar algebras over $\C$ are
considered there, the proofs of some results we need from there are
seen to carry over verbatim for general $\k$. We isolate one such fact
from \cite{Lnd} as 
a lemma below, which we will need to use repeatedly in the sequel.

\begin{Lemma}\label{zeph}
In any `exchange relation planar algebra' $P$, and for any $k \in Col$,
the space $P_k$ is linearly spanned by the images of labelled $k$-tangles
with at most $k-1$ internal boxes (all of colour $2$ and where we take $k-1 = 
0$ if
$k = 0_{\pm}$)
%(with only $2$-boxes) 
that have no `internal faces'.
\end{Lemma}

Rather than repeating the definition of an `exchange relation planar
algebra' here, it will suffice for the reader to know that if $L = L_2
= H$ and $R_0$ is any set of relations which contains $R \setminus
\{U,I\}$, then $P(L,R_0)$ is an exchange relation planar algebra.

\begin{Proposition}\label{prop:pconn} The planar algebra $P$ is a connected 
planar algebra.
\end{Proposition}

\begin{proof} We deduce from Lemma
  \ref{zeph} that $P_{0_\pm}$ is linearly spanned by its 
identity element $1_{0_\pm}$ (since $1^{0_\pm}$ is the only $0_\pm$-tangle 
`without internal faces'). Thus, $dim~P_{0_\pm} \leq 1$.

To prove the reverse inequality, we construct a functional $\lambda_\pm :
P(L)_{0_\pm} \rightarrow \k$ which is non-trivial and show that it descends to 
the quotient
$P_{0_\pm}$ of $P(L)_{0_\pm}$ by $I(R)_{0_\pm}$. The motivation for the
definition of $\lambda_\pm$ comes from the description of the partition
function of $P(H)$ given in \cite{Lnd} for $H = {\mathbb C}[G]$ and
generalised in \cite{KdyLndSnd} for a Kac algebra $H$.

The functional $\lambda_\pm$ is defined by specifying it on a basis
of $P(L)_{0_\pm}$, a typical element of which is an $L$-labelled 
$0_\pm$-tangle. Since $L_k = \emptyset$ for $k \neq 2$, each internal
box of such a tangle is necessarily a $2$-box. It will be easiest to 
illustrate the prescription in a particular example. Consider, for
instance, the $0_+$-tangle $T$ shown on the left in Figure \ref{fig:example}.
(We will identify $0_+$-tangles (resp. $0_-$-tangles) with planar networks 
with unbounded
region of colour white (resp. black) by removing the external box.)
\begin{figure}[!h]
\begin{center}
\psfrag{epa}{\huge $\epsilon(a)$}
\psfrag{eq}{\huge $=$}
\psfrag{delinphia}{\huge $\delta^{-1} \phi(a)$}
\psfrag{a1}{\huge $a_1$}
\psfrag{b1}{\huge $b_1$}
\psfrag{c1}{\huge $c_1$}
\psfrag{d1}{\huge $d_1$}
\psfrag{sa2}{\huge $Sa_2$}
\psfrag{sb2}{\huge $Sb_2$}
\psfrag{sc2}{\huge $Sc_2$}
\psfrag{sd2}{\huge $Sd_2$}
\psfrag{b}{\huge $b$}
\psfrag{a}{\huge $a$}
\psfrag{c}{\huge $c$}
\psfrag{d}{\huge $d$}
\psfrag{sa}{\huge $Sa$}
\psfrag{del}{\huge $\delta$}
\psfrag{Teq}{\huge $T = $}
\resizebox{12.0cm}{!}{\includegraphics{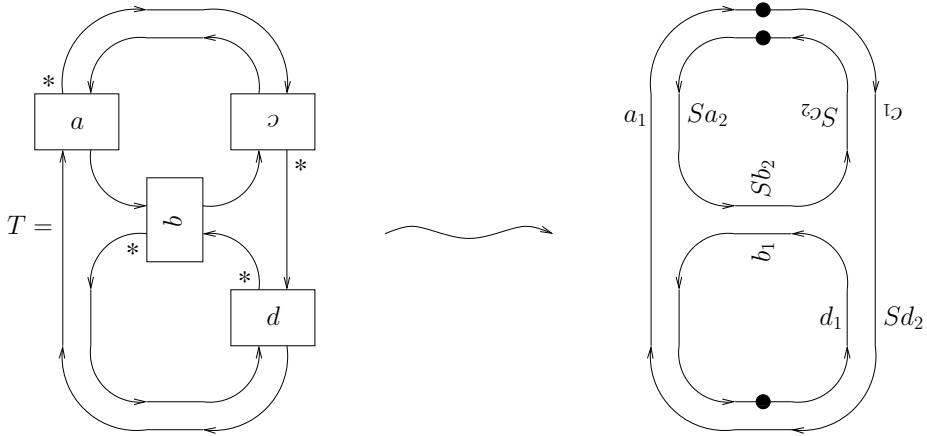}}
\end{center}
\caption{The procedure for calculating $\tau_+$}
\label{fig:example}
\end{figure} 

Removing each labelled $2$-box - with label $l$ (say) - and inserting
a symbol $l_1$ close to the strand through the $*$ of the box and a
symbol $Sl_2 = S(l_2)$ close to the other strand yields the picture
on the right of Figure \ref{fig:example}.

Now arbitrarily pick a base point on each component (loop) of the
resulting figure, read the labels on that component in the order opposite 
to that prescribed by the orientation of the loop, evaluate $\delta^{-1}\phi$
on each resulting product (the empty `product' being $1$ if the loop
has no labels) and multiply the answers.  Thus, in our example, we would obtain
\[\delta^{-1} \phi(a_1(Sd_2)c_1)
~\delta^{-1} \phi((Sc_2)(Sb_2)(Sa_2)) ~\delta^{-1} \phi(b_1d_1) 
%= \delta^{-3} \phi(a_1Sd_2c_1) \phi(Sc_2Sb_2Sa_2) \phi(b_1d_1)
~,  \]
and this element of $\k$ is what we will define as $\lambda_+(T)$. 
Note that the answer is independent of the choices of base-points
since $\phi$ is a trace. The same procedure is used to define $\lambda_-$
on $P(L)_{0_-}$.
Observe that two $0$-tangles (possibly one $0_+$ and the other $0_-$) whose
associated planar networks are isotopic on the sphere $S^2$ yield the
same element of $\k$ under the appropriate $\lambda$'s.

We now assert that $\lambda_\pm$ vanishes on $I(R)_{0_\pm}$. To see this,
it suffices to see that if two (linear combinations of) $L$-labelled 
$0_\pm$-tangles differ by an
application of any of the relations in $R$, then $\lambda_\pm$ assigns
the same value to both. For all but relation (I), this follows - and we leave 
it to the reader to verify - from
various properties of and identities in $H$ which we will list out.
\begin{itemize}
\item Relation (L) : Linearity of $(id \otimes S) \circ \Delta$.
\item Relation (M) : $\phi(1) = \delta^2$.
\item Relation (U) : $\Delta(1) = 1 \otimes 1$.
\item Relation (C) : $a_1Sa_2 = \epsilon(a).1_H$.
\item Relation (T) : $a_1\phi(Sa_2) = \phi(a).1_H$.
\item Relation (E) : $a_1 \otimes b_1 \otimes Sb_2Sa_2 = a_1 \otimes 
Sa_2(a_3b_1) \otimes S(a_4b_2)$.
\item Relation (A) : $(Sa)_1 \otimes S(Sa)_2 = Sa_2 \otimes a_1$.
\end{itemize}

The verification that $\lambda_{\pm}$ assigns the same value to two 
$L$-labelled $0_{\pm}$-tangles that differ
at a single $2$-box by an application of the relation (I) needs
a little work.
First, use isotopy on $S^2$ to move the point at infinity to a point in
the white region near the $*$ of the special $2$-boxes of both (at which they
differ). It should then be clear that the two tangles necessarily have the 
forms in Figure \ref{fig:twoforms} where $X$ is some $L$-labelled 
$2$-tangle.
\begin{figure}[!h]
\begin{center}
\psfrag{epa}{\huge $\epsilon(a)$}
\psfrag{eq}{\huge $=$}
\psfrag{delinphia}{\huge $\delta^{-1} \phi(a)$}
\psfrag{x}{\huge $X$}
\psfrag{a2b}{\huge $a_2\,b$}
\psfrag{h}{\huge $h$}
\psfrag{a}{\huge $a$}
\psfrag{sa}{\huge $Sa$}
\psfrag{t1eq}{\huge $T_1 = $}
\psfrag{t2eq}{\huge $T_2 = $}
\resizebox{8.0cm}{!}{\includegraphics{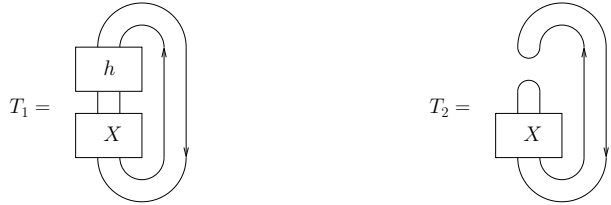}}
\end{center}
\caption{Two tangles that differ by the relation (I)}
\label{fig:twoforms}
\end{figure}
What we need to verify is that $\lambda_+(T_1) = \delta \lambda_+(T_2)$, where
the $\delta$ factor comes from the relation (I).

Let $R_0 = R \setminus \{I\}$. An application of Lemma \ref{zeph} to
$P(L,R_0)_2$ shows now that the image of
$X$ in $P_2$ may be expressed, using the relations {\it other than (I)},
as a linear combination of (images of) labelled $2$-boxes and the
tangle ${\mathcal E}^2$ shown in Figure~\ref{fig:e2}. (Tangles $\CE^3$
and $\CE^4$ are shown, in order to indicate a whole sequence of
tangles $\CE^k$.)
\begin{figure}[!h]
\begin{center}
\psfrag{epa}{\huge $\epsilon(a)$}
\psfrag{eq}{\huge $=$}
\psfrag{delinphia}{\huge $\delta^{-1} \phi(a)$}
\psfrag{a1}{\huge $a_1$}
\psfrag{a2b}{\huge $a_2\,b$}
\psfrag{b}{\huge $b$}
\psfrag{a}{\huge $a$}
\psfrag{sa}{\huge $Sa$}
\psfrag{del}{\huge $\delta$}
\resizebox{6cm}{!}{\includegraphics{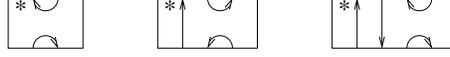}}
\end{center}
\caption{The tangles ${\mathcal E}^2,{\mathcal E}^3$ and ${\mathcal E}^4$}
\label{fig:e2}
\end{figure}
Since we have already verified that $\lambda_\pm$ is invariant under
application of any of the relations other than (I), we reduce immediately
to the case that $X$ itself is either a labelled $2$-box or the tangle
${\mathcal E}^2$. We treat these cases one by one.

If $X$ is a $2$-box labelled by $a \in H$, then the procedure for calculating
$\lambda_+$ yields $\lambda_+(T_1) = \delta^{-1} \phi(h_1a_1) \delta^{-1}
\phi(Sh_2Sa_2) = \delta^2 \epsilon(a)$ while $\lambda_+(T_2) = \delta^{-1}
\phi(Sa_2a_1) = \delta \epsilon(a)$. 
On the other hand, if $X = {\mathcal E}^2$, then $\lambda_+(T_1)$ may be
computed to be $\delta^{-1} \phi(h_1Sh_2) = \delta^{3}$ while
$\lambda_+(T_2)$ is computed to be $\delta^{-2} \phi(1)\phi(1) = \delta^2$.
In either case, we see that $\lambda_+(T_1) = \delta \lambda_+(T_2)$.
This completes verification of invariance under the relation (I).

Thus $\lambda_\pm$ descend to give maps from $P_{0_\pm} \rightarrow \k$
that are clearly (since $\lambda_\pm(1_{0_\pm}) = 1$) surjective, thereby
establishing that $dim ~P_{0_\pm} \geq 1$ and concluding the proof.
\end{proof}

Before stating our next result, note that since for any planar algebra $P$, 
each $P_k$ is a $\k$-algebra with identity $1_k$, it
follows that if $P$ is connected, there are canonical 
identifications
$P_{0_\pm} = \k$ (under which $1_{0_\pm}$ is identified with $1$). 
Consequently, if $N$ is any planar network  then its 
partition function 
$Z^{P}_N$ takes values in $\k$. 
%(This identification is, in fact,
%implicit in the relation (M).)

\begin{Corollary}
For any labelled planar network $N$ in $P(L)_{0_\pm}$, its value (as given
by the partition function of $P$) is  $\lambda_\pm(N)$.
\end{Corollary}

\begin{proof}
By Lemma \ref{zeph}, any $L$-labelled $0_\pm$-tangle may be written,
using all the relations in $R$, as a multiple of $1^{0_\pm}$ - the
only $0_\pm$-tangle without internal faces. Since $\lambda_\pm$ is 
invariant under the relation $R$, it suffices to verify that for the
tangle $1^{0_\pm}$, both the partition function and $\lambda_\pm$ 
assign the same value. This is true since both give the value $1$ to
$1^{0_\pm}$, establishing the corollary.
\end{proof}

\begin{figure}[!h]
\begin{center}
\psfrag{epa}{\huge $\epsilon(a)$}
\psfrag{eq}{\huge $=$}
\psfrag{delinphia}{\huge $\delta^{-1} \phi(a)$}
\psfrag{a1}{\huge $a_1$}
\psfrag{a2b}{\huge $a_2\,b$}
\psfrag{b}{\huge $b$}
\psfrag{a}{\huge $a$}
\psfrag{sa}{\huge $Sa$}
\psfrag{del}{\huge $\delta$}
\resizebox{1.5cm}{!}{\includegraphics{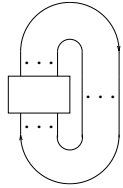}}
\end{center}
\caption{The trace tangle}
\label{fig:trtgl}
\end{figure}

Since we have now verified that $P_{0_\pm}$ is $1$-dimensional, and in 
particular that its identity 
$1_{0_\pm}$  is non-zero, the following equation for $x \in P_k$ where $k \geq 
1$,
\[ Z^{~\,P}_{Tr(k)}(x) = \delta^k \tau_k(x) 1_{0_+}, \]
uniquely
defines a tracial linear functional $\tau_k$ on $P_k$  
which will
be referred to as the normalised picture trace. Here $Tr(k)$ denotes 
the $0_+$-tangle with a single internal $k$-box that is shown in Figure
\ref{fig:trtgl}.
Note that $\tau_k(1_k) = 1$.

Let $\pi$ denote the natural map from $H$ to $P_2$ which takes $a \in H$
to the image of the $2$-box labelled by $a$. We then have the following lemma
by which we will henceforth identify $H$ with $P_2$.

\begin{Lemma}\label{lemma:pilemma}
The map $\pi$ is a unital algebra isomorphism from $H$ to $P_2$.
\end{Lemma}

\begin{proof}
From the results of \cite{Lnd} it follows that each element of $P_2$
is a linear combination of labelled $2$-boxes or equivalently that $\pi$
is surjective. The relations (E) and (C) may be seen to imply that
$\pi$ preserves multiplication while (U) ensures that $\pi$ is unital.

To show that $\pi$ is injective, observe that relations (T) and (M)
imply that $\tau_2 \circ \pi = \delta^{-2}\phi$.
Since $\delta^{-2}\phi$ is a non-degenerate trace on $H$, if  $\{e_i:i\in I\}$ 
is a basis of $H$ there is a basis $\{e^i: i \in
I\}$ of $H$ that is dual to this basis in the sense that 
$\delta^{-2}\phi(e_ie^j) = \delta_i^j$ - where, of course, $\delta_i^j$ is the 
Kronecker delta.
Since $\pi$ preserves multiplication, it follows that 
$\tau_2(\pi(e_i)\pi(e^j)) = \delta_i^j$. Thus $\{\pi(e_i) : i \in I\}$
is linearly independent and so $\pi$ is a unital algebra isomorphism.
\end{proof}

We omit the proof of the following corollary which follows immediately from
Lemma \ref{lemma:pilemma} and the relation (T).

\begin{Corollary}
The planar algebra $P$ is irreducible.\qed
\end{Corollary}

\begin{figure}[!h]
\begin{center}
\psfrag{epa}{\huge $\epsilon(a)$}
\psfrag{eq}{\huge $=$}
\psfrag{delinphia}{\huge $\delta^{-1} \phi(a)$}
\psfrag{a1}{\huge $a_1$}
\psfrag{a2b}{\huge $a_2\,b$}
\psfrag{1}{\huge $1$}
\psfrag{2}{\huge $2$}
\psfrag{3}{\huge $3$}
\psfrag{4}{\huge $4$}
\psfrag{5}{\huge $5$}
\psfrag{a}{\huge $a$}
\psfrag{sa}{\huge $Sa$}
\psfrag{del}{\huge $\delta$}
\resizebox{10.0cm}{!}{\includegraphics{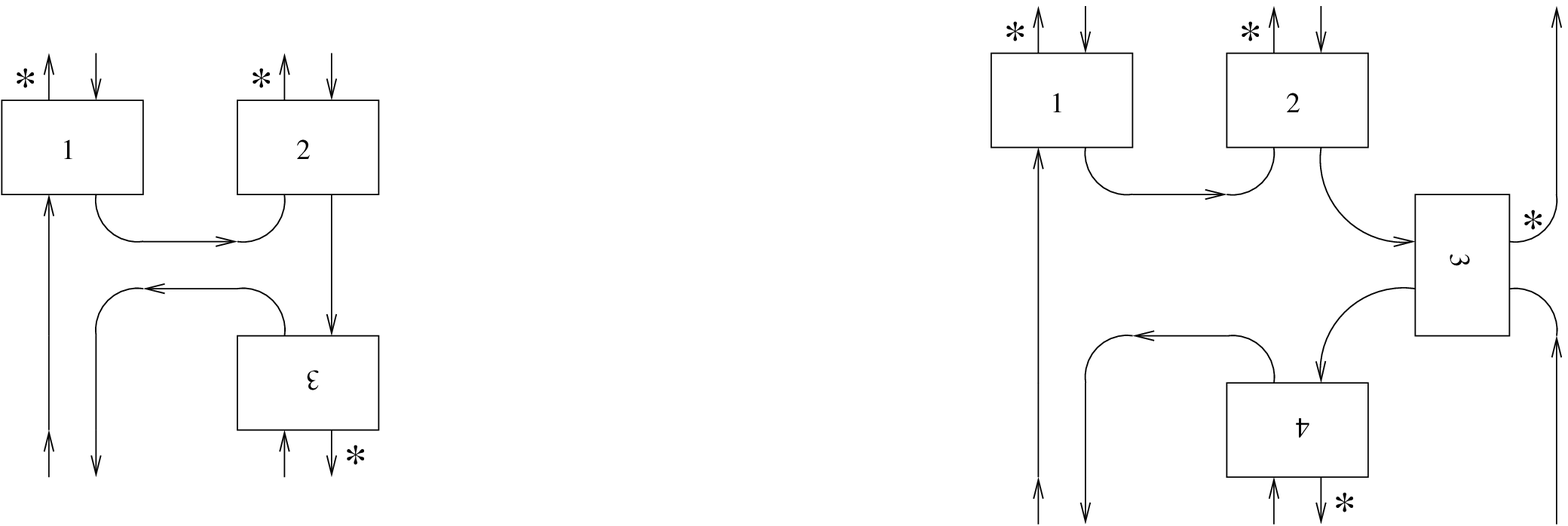}}
\end{center}
\caption{The tangles $X_4$ and $X_5$}
\label{fig:xktangle}
\end{figure}

\begin{figure}[!h]
\begin{center}
\psfrag{epa}{\huge $\epsilon(a)$}
\psfrag{eq}{\huge $=$}
\psfrag{delinphia}{\huge $\delta^{-1} \phi(a)$}
\psfrag{a1}{\huge $a_1$}
\psfrag{a2b}{\huge $a_2\,b$}
\psfrag{1}{\huge $1$}
\psfrag{2}{\huge $2$}
\psfrag{3}{\huge $3$}
\psfrag{4}{\huge $4$}
\psfrag{5}{\huge $5$}
\psfrag{a}{\huge $a$}
\psfrag{sa}{\huge $Sa$}
\psfrag{del}{\huge $\delta$}
\resizebox{10.0cm}{!}{\includegraphics{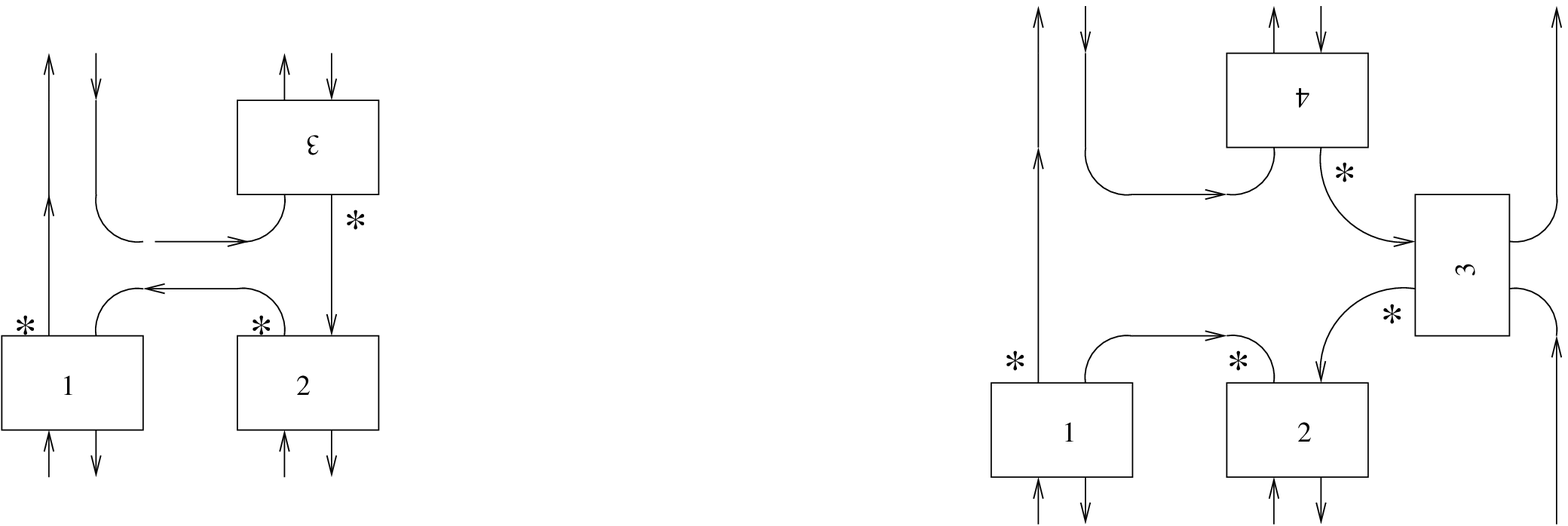}}
\end{center}
\caption{The tangles $X_4^*$ and $X_5^*$}
\label{fig:xksttangle}
\end{figure}

\begin{figure}[!h]
\begin{center}
\psfrag{xk}{\huge $X_k$}
\psfrag{kp1}{\huge $k$}
\psfrag{xkp1}{\huge $X_{k+1}$}
\resizebox{3.5cm}{!}{\includegraphics{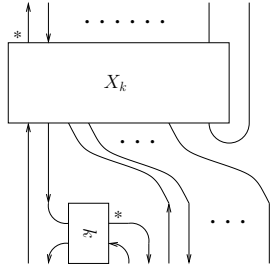}}
\end{center}
\caption{Inductive definition of $X_{k+1}$}
\label{fig:xkind}
\end{figure}

We will be interested in describing a basis of $P_k$ in terms of a
basis of $H$. For this, the $k$-tangles $X_k$ and their adjoint tangles
$X_k^*$ (each with $k-1$ internal $2$-boxes), illustrated in Figures 
\ref{fig:xktangle} and \ref{fig:xksttangle} for $k=4$ and $k=5$, will turn out 
to be relevant.
Note that the tangles $X_k$ may be defined inductively as in Figure 
\ref{fig:xkind}.
Close relatives of the tangles $X_k$ occur in \cite{LndSnd}.

Before describing a basis of $P_k$, we would like to introduce the notation 
${\mathcal T}(k)$ (resp. ${\mathcal T}_\leq(k)$) for the set of $k$-tangles
without internal faces and exactly $k-1$ (resp. at most $k-1$) internal
boxes all of colour $2$. Note that $X_k, X_k^* \in {\mathcal T}(k)$.

We will next prove the following lemma. See the appendix for a generalisation.
\begin{Lemma}\label{lemma:surj}
For each $k \geq 3$, the
map $Z^P_{X_k} : P_2^{\otimes (k-1)} \rightarrow P_k$ is surjective.
\end{Lemma}

\begin{proof} It is easy to verify that any element of ${\mathcal T}(3)$ is
obtained from any other by rotating the external and internal
boxes. (cf. eq. (3.1) in \cite{KdyLndSnd}.)
It then follows from relations (E) and (A) (and (L)) that all these tangles 
have the
same range in $P_3$.
Further, inspection shows that any element of ${\mathcal T}_\leq(3)$ may be 
obtained
from one of  the elements of ${\mathcal T}(3)$ by substituting $1^2$ or
${\mathcal E}^2$ into some of its internal boxes.
Together with Lemma \ref{zeph},
this implies that for any $X \in {\mathcal T}(3)$, the map $Z^P_{X}$ is 
surjective.
In particular, $Z^P_{X_3}$ is surjective.

The tangles $X$ and $W$ of Figure \ref{fig:p2e2p2} are in ${\mathcal T}(3)$ 
and the
surjectivity of $Z^P_{X}$ implies that
\begin{figure}[!h]
\begin{center}
\psfrag{xk}{\huge $X_k$}
\psfrag{xeq}{\huge $X =$}
\psfrag{weq}{\huge $W =$}
\psfrag{feq}{\huge $F =$}
\psfrag{leq}{\huge $L =$}
\psfrag{kp1}{\huge $k+1$}
\psfrag{xkp1}{\huge $X_{k+1}$}
\resizebox{12.0cm}{!}{\includegraphics{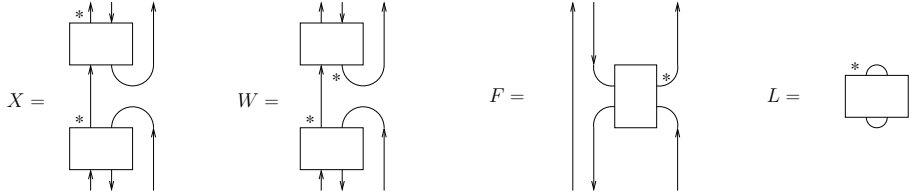}}
\end{center}
\caption{The tangles $X$, $W$, $F$ and $L$}
\label{fig:p2e2p2}
\end{figure}
$P_3 = P_2e_2P_2$ and consequently (see 
Lemma 5.7 of \cite{KdyLndSnd} for the proof) that $P_{k+1} = 
P_2e_2e_3\cdots e_{k}P_{k}$ for any $k \geq 3$. 
(We recall that $e_k =
\delta^{-1} Z_{\CE^{k+1}}(1)$, with $\CE^{k+1}$ as indicated by Figure
\ref{fig:e2}.)

This is equivalent to the statement that the $k+1$ tangle $T$ with two 
internal boxes of colours $2$ and
$k$ shown in Figure \ref{fig:k2ml} has $Z_T^P$ surjective. 
\begin{figure}[!h]
\begin{center}
\psfrag{1}{\huge $1$}
\psfrag{2}{\huge $2$}
\psfrag{xk}{\huge $X_k$}
\psfrag{kp1}{\huge $k+1$}
\psfrag{xkp1}{\huge $X_{k+1}$}
\resizebox{3.5cm}{!}{\includegraphics{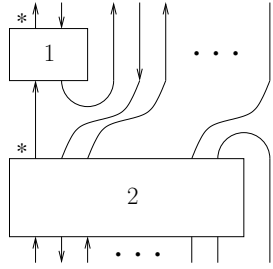}}
\end{center}
\caption{The $k+1$-tangle $T$}
\label{fig:k2ml}
\end{figure}
However, a little thought shows that $R_{k+1}^{-2} \circ X_{k+1} = 
T \circ_{D_2} (R_{k}^{-1} \circ X_{k})$ where $R_k$ is the $k$-rotation
tangle.
Since $Z_{R_k}$ is an isomorphism for each $k$, it follows that surjectivity 
of  $Z^P_{X_{k}}$
implies that of $Z^P_{X_{k+1}}$ and concludes the proof, by induction on $k$.
\end{proof}

In order to state our next result, we
fix a pair of bases $\{e_i:i\in I\}$ and 
$\{e^i:i\in I\}$ of $H$ that are dual with respect to the trace 
$\delta^{-2}\phi$, as in Lemma \ref{lemma:pilemma}. We then have:

\begin{Lemma}\label{lemma:eiej}For $k \geq 2$, and for each ${\bf i} =
(i_1,\cdots,i_{k-1}) \in I^{k-1}$, if we set $e_{{\bf i}} =
Z_{X_k}(e_{i_1}\otimes \cdots \otimes  e_{i_{k-1}})$ and
$e^{{\bf i}} = Z_{X_k^*}(e^{i_1}\otimes  
\cdots \otimes  e^{i_{k-1}})$, then $\{e_{\bf i}:{\bf i} \in 
I^{k-1}\}$ 
and 
$\{e^{\bf i}:{\bf i} \in I^{k-1}\}$ are a pair of bases of $P_k$
dual with respect to the trace $\tau_k$, which is a non-degenerate trace.
In particular, $dim~P_k = n^{k-1}$.
% and $\pi_k$ is a $\k$-vector space isomorphism.
\end{Lemma}

\begin{proof}
The case $k=2$ of this lemma is contained in Lemma \ref{lemma:pilemma}.
For a general $k$, Lemma \ref{lemma:surj} shows that
$P_k$ 
is linearly spanned by (the images of) $H$-labelled
$k$-tangles $Z_{X_k}(a(1)\otimes \cdots \otimes  a(k-1)
\otimes 1)$ where 
$(a(1),\cdots,a(k-1)) \in H^{k-1}$.  
This establishes the inequality $dim ~P_k \leq n^{k-1}$.

We will next show that for $k=4$, if ${\bf i} = (i_1,i_2,i_3)$ and ${\bf j} =
(j_1,j_2,j_3)$, then $\tau_4(e_{{\bf i}}e^{{\bf
    j}}) = \delta^{\bf j}_{\bf i}$; the proof of the general case is similar.
Notice that $\tau_4(e_{{\bf i}}e^{{\bf
    j}}) = \delta^{-4} \times$(the value of the labelled planar network $N$ in 
Figure \ref{fig:t4eiej}):

\begin{figure}[!h]
\begin{center}
\psfrag{epa}{\huge $\epsilon(a)$}
\psfrag{eq}{\huge $=$}
\psfrag{delinphia}{\huge $\delta^{-1} \phi(a)$}
\psfrag{a1}{\huge $a_1$}
\psfrag{a2b}{\huge $a_2\,b$}
\psfrag{1}{\huge $e_{i_1}$}
\psfrag{2}{\huge $e_{i_2}$}
\psfrag{3}{\huge $e_{i_3}$}
\psfrag{4}{\huge $e^{j_3}$}
\psfrag{5}{\huge $e^{j_2}$}
\psfrag{6}{\huge $e^{j_1}$}
\psfrag{a}{\huge $a$}
\psfrag{sa}{\huge $Sa$}
\psfrag{del}{\huge $\delta$}
\resizebox{5.0cm}{!}{\includegraphics{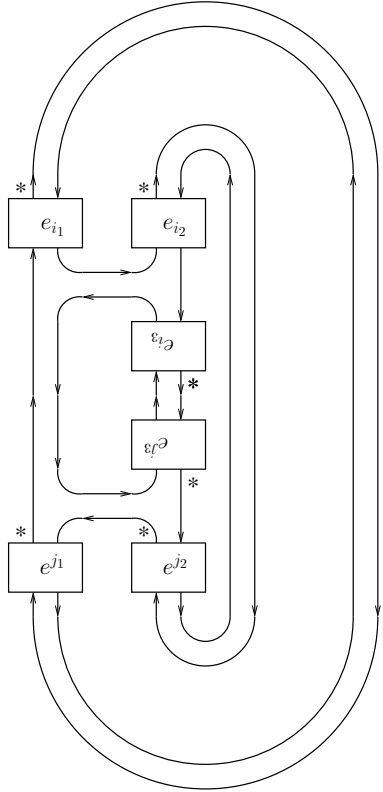}}
\end{center}
\caption{The planar network $N$ of Lemma \ref{lemma:eiej}}
\label{fig:t4eiej}
\end{figure}
Now use relation (T) and the hypothesis that $\delta^{-2}\phi(e_ie^j) =
\delta^{-2}\phi(e^je_i) = \delta_i^j$ to conclude that 
the relation shown in Figure \ref{fig:eiej} holds in $P$.
\begin{figure}[!h]
\begin{center}
\psfrag{epa}{\huge $\epsilon(a)$}
\psfrag{eq}{\huge $=$}
\psfrag{delinv}{\huge $\delta^{-1}$}
\psfrag{ejei}{\huge $e^je_i$}
\psfrag{t2ejei}{\huge $\delta^{-2}\phi(e^je_i)$}
\psfrag{1}{\huge $e_{i_1}$}
\psfrag{2}{\huge $e_{i_2}$}
\psfrag{3}{\huge $e_{i_3}$}
\psfrag{4}{\huge $e^{j_3}$}
\psfrag{5}{\huge $e^{j_2}$}
\psfrag{6}{\huge $e^{j_1}$}
\psfrag{a}{\huge $a$}
\psfrag{sa}{\huge $Sa$}
\psfrag{delij}{\huge $\delta^j_i$}
\resizebox{8.0cm}{!}{\includegraphics{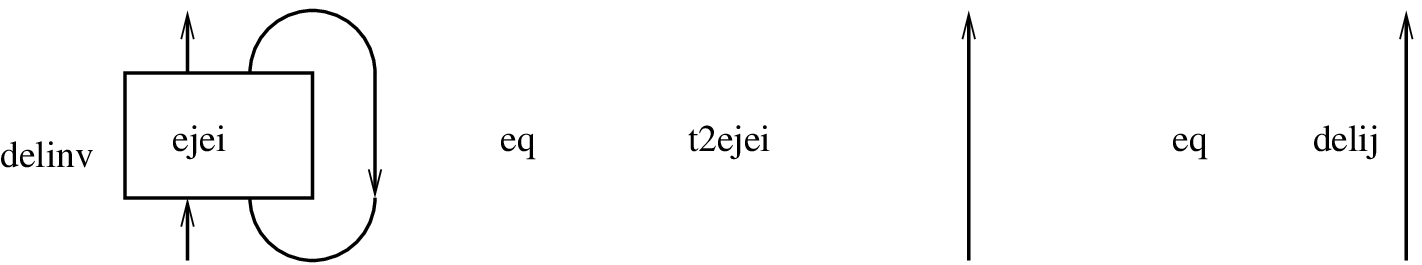}}
\end{center}
\caption{A relation that holds in $P$}
\label{fig:eiej}
\end{figure}
Apply this repeatedly with $(i,j) = 
(i_3,j_3), (i_2,j_2), (i_1,j_1)$ to deduce - after $4$ applications
of relation (M) that
$\tau_4(e_{{\bf i}}e^{{\bf j}}) = \delta^{{\bf j}}_{{\bf
    i}},$
thereby establishing the desired equation.

It now follows that $\{e_{\bf i}:{\bf i} \in I^{k-1}\}$ is a linearly
independent set in $P_k$, that $dim ~P_k = n^{k-1}$, that $\{e_{\bf i}:{\bf i} 
\in I^{k-1}\}$ is a basis of $P_k$ with $\{e^{\bf i}:{\bf i} \in I^{k-1}\}$
being a dual basis, and finally that $\tau_k$ is a non-degenerate trace
on $P_k$.
\end{proof}

Recall that a planar algebra $P$ is said to be spherical if the 
partition function for
planar networks is an invariant of isotopy on $S^2$.
We now make
the following simple observation.

\begin{Lemma} 
If $P$ is any connected and  irreducible planar algebra with modulus $\delta$, 
then $P$ is spherical.
\end{Lemma}

\begin{proof}
Since the only difference between viewing a network as being embedded
in the plane or on the sphere is how it is positioned with respect to
the point at infinity, it is seen after a little thought that a
connected planar
algebra is spherical if and only if $Z^P_{T_L} = Z^P_{T_R}$ where $T_L$
(resp., $T_R$) is the $0_-$-tangle (resp., $0_+$-tangle) shown in Figure \ref{fig:trtgls},
where both $Z^{P}_{T_L}$ and $Z^{P}_{T_R}$ are regarded as
linear functionals on $P_1$:
\begin{figure}[!h]\label{fig:trtgls}
\begin{center}
\psfrag{epa}{\huge $\epsilon(a)$}
\psfrag{eq}{\huge $=$}
\psfrag{delinphia}{\huge $\delta^{-1} \phi(a)$}
\psfrag{a1}{\huge $a_1$}
\psfrag{a2b}{\huge $a_2\,b$}
\psfrag{eqtr}{\huge $T_R = $}
\psfrag{tleq}{\huge $T_L = $}
\psfrag{sa}{\huge $Sa$}
\psfrag{del}{\huge $\delta$}
\resizebox{8cm}{!}{\includegraphics{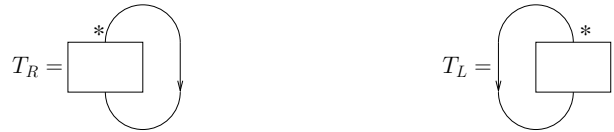}}
\end{center}
\caption{The tangles $T_R$ and $T_L$}
%\label{fig:trtgls}
\end{figure}

However, for irreducible $P$, the space $P_1$ is 1-dimensional
and is consequently spanned by its identity element $1_1$, and 
relation (M) says precisely that 
\[Z_{T_L}^P(1_1)  = \delta~ = Z_{T_R}^P(1_1);\]
and since two linear functionals which agree on a basis must be
identical, we see that $P$ is indeed spherical.
\end{proof}  

We summarise the facts we have proved about $P$ in the following theorem.
The term `non-degenerate planar algebra' is used for a connected planar
algebra for which the picture traces $\tau_k$ are all non-degenerate.
Recall that a planar algebra $P$ with non-zero modulus is said to be of depth 
two if $P_3
= P_2e_2P_2$, where $e_2$ is defined as in the proof of Lemma \ref{lemma:surj} 
(or equivalently, if $Z^P_X$ is 
surjective where $X$ is the tangle of 
Figure \ref{fig:p2e2p2}).

\begin{Theorem}
Let $H$ be a semisimple and
cosemisimple Hopf algebra $H$ of dimension $n$.
The planar algebra $P = P(H,\delta)$ associated  to $H$ is a connected, 
irreducible, 
spherical, non-degenerate planar algebra with modulus $\delta$ and of depth 
two. Further, 
$dim ~P_k = n^{k-1}$ for all $k \geq 1$.\qed 
\end{Theorem}

It should be clear that an isomorphism of semisimple and cosemisimple
Hopf algebras naturally yields an isomorphism of the corresponding
planar algebras (with the same choice of $\delta$).

\section{From planar algebras to Hopf algebras}

In this section, we wish to invert the procedure of \S 3 to get
a semisimple and cosemisimple Hopf algebra from a connected, irreducible and
non-degenerate planar algebra of depth two and non-zero modulus.
So fix such a planar algebra $P$ with modulus $\delta$.

If $P$ is a `subfactor planar algebra', there is a detailed description
in \cite{DasKdy} of the construction of a Kac algebra from $P$. Essentially
the same proof works in our situation to get a Hopf algebra from $P$,
and so we will only indicate the changes to be made for that proof to
work here. These changes are summarised in the following lemmas.

\begin{Lemma}
For the tangle $W$ shown in Figure \ref{fig:p2e2p2}, the map $Z_W^P : P_2
\otimes P_2 \rightarrow P_3$ is an isomorphism.
\end{Lemma}

\begin{proof}
The depth two assumption on $P$ is equivalent to the surjectivity of
$Z_X^P$ for the tangle $X$ of Figure \ref{fig:p2e2p2} or to that of
$Z_W^P$ - since $W = X \circ_{D_1} R_2$ and all rotation tangles
give isomorphisms of the spaces they naturally act on. Thus $Z_W^P$
is surjective.

As for injectivity, use the non-degeneracy of the picture trace to
choose a pair of bases $\{e_i : i \in I\}$ and $\{e^i : i \in I\}$
of $P_2$ that are dual with respect to $\tau_2$.
Then, the proof of Lemma \ref{lemma:eiej} goes through to show that
$\{ Z_{X_3}(e_{i_1} \otimes e_{i_2}): i_1,i_2 \in I\}$ forms a linearly 
independent set in $P_3$. Thus $dim~P_3 = (dim~P_2)^2$ and so
$Z_W^P$ is an isomorphism. 
\end{proof}

The proof of the main theorem of \cite{DasKdy} goes through in this
context to imply that $H = P_2$ is 
a Hopf algebra with its usual (from the planar algebra) multiplicative 
structure
and the comultiplication, counit and antipode defined by $\Delta(a) 
= (Z^P_W)^{-1}Z_F^P(a)$, $\epsilon(a) = \delta^{-1}Z_L^P(a)$ (where the
tangles $F$ and $L$ are shown in Figure \ref{fig:p2e2p2}) and $S(a) = 
Z_{R_2}^P(a)$. We now have:
%The next lemma assures us that the $H$ obtained as above is
%indeed semisimple and cosemisimple.

\begin{Lemma}
With the foregoing notations, the Hopf algebra $H$ is semisimple and
cosemisimple and its dimension $n$ is related to the modulus $\delta$
of $P$ by $\delta^2 = n$.
\end{Lemma}

\begin{proof}
Recall that a Hopf algebra $H$ is semisimple if there exists a one-sided
integral $h \in H$ with $\epsilon(h) \neq 0$.
For $H = P_2$, define $h \in H$ and $\phi \in H^*$ so that
the (I) and (T) relations hold (for the $\phi$ this needs irreducibility
of $P$).

A pleasant exercise with the relations in Figures \ref{fig:LnrMdl} - 
\ref{fig:XchNtp}
then shows that $ha = \epsilon(a)h$, $a_1\phi(a_2) = \phi(a).1$,
$\epsilon(h) = \delta^2$, $\phi(1) = \delta^2$ and $\phi(h) = \delta^2$ 
proving that both $H$ and $H^*$ are semisimple.
However, in a semisimple and cosemisimple Hopf algebra, there are
choices of $h$ and $\phi$ for which $\epsilon(h) = \phi(1) = \phi(h) =n$;
since the space of integrals is 1-dimensional, it follows
that $\delta^2 = n$. 
\end{proof}

\begin{Proposition}\label{isom}
The association $H \mapsto P(H,\delta)$ defines a bijective
correspondence between isomorphism classes of
semisimple and cosemisimple Hopf algebras (over $\k$) with $dim ~H =
\delta^2 \in \k$, on the one
hand, and isomorphism classes of 
connected, irreducible, non-degenerate planar algebras (over $\k$)
with modulus $\delta$ and of depth two.
\end{Proposition}

\begin{proof} It is easy to see that 
\[ H_1 \cong H_2 \Rightarrow P(H_1,\delta) \cong P(H_2,\delta)~.\]

In the other direction, suppose $P$ is a connected, irreducible,
non-degenerate planar algebra (over $\k$) with modulus $\delta$ and of
depth two. Let $H$ be the semisimple and cosemisimple Hopf algebra
constructed as above. We wish to prove first that $P \cong P(H,\delta)$.

Since $P_2 = H$, there is a planar algebra homomorphism of $\pi: P(L)
\rightarrow P$ -
where $L = L_2 = H$. The depth two assumption says $P_3 = P_2
e_2P_2$, which implies (as already observed in the proof of Lemma
\ref{lemma:surj}) that $P_{k+1} = P_2e_2e_3\cdots e_kP_k ~\forall k
\geq 3$, and hence (by induction) that $P$ is generated, as a planar
algebra, by $P_2$; and in particular, the map $\pi$ is surjective.

Next, it is easy to see that all the relations defining $P(H,\delta)$
are satisfied in $P$, and hence $\pi$ descends to a surjective planar
algebra homomorphism of $P(H,\delta)$ to $P$. In particular, $dim ~P_k
\leq dim ~P(H,\delta)_k = (dim ~H)^{k-1}$. On the other hand, the
proof of Lemma \ref{lemma:eiej} shows, even in this case, that if
$\{e_i:i \in I\}$ and $\{e_j:j \in I\}$ are a pair of bases of $P_2$ which
are dual with respect to $\tau_2$, then $\{e_{\bf i}: {\bf i} \in I^{k-1}\}$
and $\{e_{\bf j}: {\bf j} \in I^{k-1}\}$ (as defined in Lemma
\ref{lemma:eiej}) are linearly independent in $P_k$ and that hence, also
$dim ~P_k \geq  (dim ~H)^{k-1} = dim ~P(H,\delta)_k$. This shows that
indeed $P(H,\delta) \cong P$.

To complete the proof, note that 
\[P(H_1,\delta) \stackrel{\psi}{\cong} P(H_2,\delta) \Rightarrow H_1
\stackrel{\psi_2}{\cong} H_2 .\] 
%with the isomorphism being given by
%\[H_1 = (P(H_1,\delta)_2 \stackrel{\psi_2}{\rightarrow}
%(P(H_2,\delta)_2 = H_2.\]
\end{proof}

\section{Duality between $P(H,\delta)$ and $P(H^*,\delta)$}
We will next explicate a duality between the planar algebras associated
to $H$ and to $H^*$. 
Recall from \cite{KdySnd} that there is an `operation
on planar tangles' denoted by `-'.
This is defined by (i) the map $k \mapsto k^-$ that toggles $0_\pm$ and fixes
the other colours and (ii) the
map $T \mapsto T^-$ that moves the $*$  back (anticlockwise)
by one on all boxes and inverts shading.
If $P$ is a planar algebra, the planar algebra $^-P$ is defined by 
setting $^-P_k = P_{k^-}$
and $Z^{^-P}_T = Z^P_{T^-}$ for each tangle $T$.
By $\#$, we will denote the inverse operation (which moves all $*$s forward by
one and inverts shading).

We will also need to recall the Fourier transform map for $H$. This is the map 
$F : H \rightarrow H^*$ defined by $F(a) = \delta^{-1} \phi_1(a) \phi_2$. We
use $F$ to also denote the Fourier transform map of $H^*$, the argument of $F$
making it clear which one is meant. Similarly, we use $S$ to also denote the 
antipode of $H^*$. The properties of $F$ that we will use are
that $F^2 = S$, $FS = SF$ and $F(SF) = id = (SF)F$.

The result that we wish to prove is the following theorem.

\begin{Theorem}\label{theo:pliso}
The map $^-P(H,\delta)_2 = H \rightarrow H^* =
P(H^*,\delta)_2$ defined by $a \mapsto SF(a)$ extends to a
planar algebra isomorphism from $^-P(H,\delta)$ to $P(H^*,\delta)$.
\end{Theorem}

\begin{proof} 
Observe first that if $P$ is a planar algebra presented with generators $L =
\coprod_{k \in Col} L_k$
and relations $R$, then the planar algebra $^-P$ is presented by the
label set given by $^-L_k = L_{k^-}$, and relations $^-R$ given as follows.
Consider a typical relation in $R$. It is given as a linear combination
of $L$-labelled tangles all of a fixed colour.
Applying $\#$ to each of these (leaving the labels unchanged but regarded
as elements of $^-L$) gives a linear combination of $^-L$-labelled tangles
which is the typical relation of $^-R$. 
This is an easy consequence of the definitions in \cite{KdySnd}.

In particular, $^-P(H,\delta)$ is presented with generators $^-L$,
where $^-L_2 = H$ and $^-L_k = \emptyset$ for $k \neq 2$, and relations
given by the $\#$ of the relations in Figures \ref{fig:LnrMdl} - 
\ref{fig:XchNtp}. For instance, the relations corresponding to those in Figure 
\ref{fig:XchNtp}
are given by those in Figure \ref{fig:XchNtpnew}.

\begin{figure}[!h]
\begin{center}
\psfrag{epa}{\huge $\epsilon(a)$}
\psfrag{eq}{\huge $=$}
\psfrag{delinphia}{\huge $\delta^{-1} \phi(a)$}
\psfrag{a1}{\huge $a_1$}
\psfrag{a2b}{\huge $a_2\,b$}
\psfrag{b}{\huge $b$}
\psfrag{a}{\huge $a$}
\psfrag{sa}{\huge $Sa$}
\psfrag{del}{\huge $\delta$}
\resizebox{12.0cm}{!}{\includegraphics{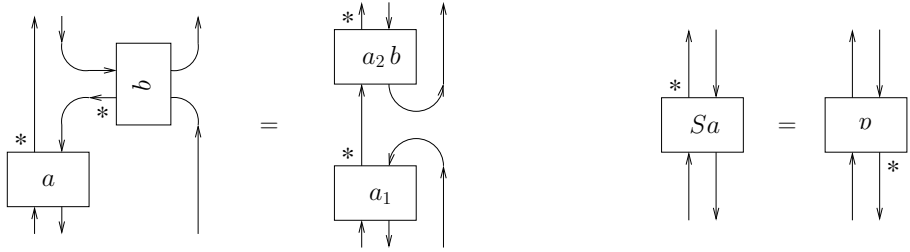}}
\end{center}
\caption{The E(xchange) and A(ntipode) relations in $^-P(H,\delta)$}
\label{fig:XchNtpnew}
\end{figure}

The universal property of the planar algebra $P(^-L)$ implies that
there is a planar algebra map (i.e., a map of vector spaces for each $k \in 
Col$ that intertwines the tangle actions) from $P(^-L)$ to $P(H^*,\delta)$
that takes the $2$-box labelled by $a$ in $P(^-L)$ to the (image of the) one 
labelled by $SF(a)$ in $P(H^*,\delta)$.
Since $P(H^*,\delta)$ is generated by its $2$-boxes as a planar algebra,
this map is surjective.
If we now verify that all relations in $^-R$ go to $0$ under this map,
it will induce a surjective planar algebra map from the quotient 
$^-P(H,\delta)$ to $P(H^*,\delta)$.
A comparison of the dimensions will then show that this is a planar algebra
isomorphism and conclude the proof.

It remains to verify that for each relation in $^-P(H,\delta)$, the relation
obtained by substituting $SF(a)$ in each box labelled by $a$ gives
a valid relation in $P(H^*,\delta)$.
As in the proof of Proposition \ref{prop:pconn}, we will leave all the
easier verifications to the reader indicating only the relevant properties and 
the appropriate relations in $P(H^*,\delta)$ used.
\begin{itemize}
\item Relation (L) : Linearity of $SF$ and relation (L).
\item Relation (M) : The equality of the choice of $\delta$ for $P(H,\delta)$
and $P(H^*,\delta)$ and relation (M).
\item Relation (U) : $SF(1_H) = \delta^{-1} \phi$ and relation (I).
\item Relation (I) : $SF(h) = \delta \epsilon$ and relation (U).
\item Relation (C) : $\delta^{-1} (SF(a))(h) = \epsilon(a)$ and relation (T).
\item Relation (T) : $(F(a))(1) = \delta^{-1}\phi(a)$ and relations
  (A) and (C).
\item Relation (A) : $SFS = F$ and relation (A).
\end{itemize}

For relation (E), it follows from the figure on the left in
Figure \ref{fig:XchNtpnew} that the relation that requires to be
verified in $P(H^*,\delta)$ is the one in Figure \ref{fig:needs},
for each $a,b \in H$.
\begin{figure}[!h]
\begin{center}
\psfrag{epa}{\huge $\epsilon(a)$}
\psfrag{eq}{\huge $=$}
\psfrag{delinphia}{\huge $\delta^{-1} \phi(a)$}
\psfrag{a1}{\Large $SF(a_1)$}
\psfrag{a2b}{\Large $SF(a_2\,b)$}
\psfrag{b}{\Large $SF(b)$}
\psfrag{a}{\Large $SF(a)$}
\psfrag{sa}{\huge $Sa$}
\psfrag{del}{\huge $\delta$}
\resizebox{7.0cm}{!}{\includegraphics{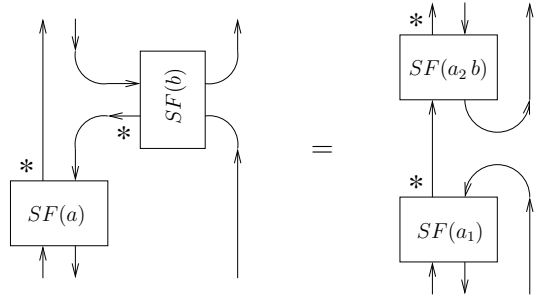}}
\end{center}
\caption{Relation to be verified in $^-P(H^*,\delta)$}
\label{fig:needs}
\end{figure}

Note that $SF(a) = \delta^{-1} \phi_1(a)S\phi_2$ and $SF(b) = \delta^{-1} 
\widetilde{\phi}_1(a)S\widetilde{\phi}_2$, where $\widetilde{\phi}$ is another 
copy of $\phi$. Then the labelled tangle
on the left of Figure \ref{fig:needs} equals (by application of the (A) and (E)
relations in $P(H^*,\delta)$) either of the labelled tangles in Figure 
\ref{fig:needs2}.
\begin{figure}[!h]
\begin{center}
\psfrag{epa}{\huge $\epsilon(a)$}
\psfrag{eq}{\huge $=$}
\psfrag{delinphia}{\huge $\delta^{-1} \phi(a)$}
\psfrag{a1}{\huge $\widetilde{\phi}_3S\phi_2$}
\psfrag{a2b}{\huge $S\widetilde{\phi}_2$}
\psfrag{b}{\huge $\widetilde{\phi}_2$}
\psfrag{a}{\huge $\phi_2$}
\psfrag{u}{\huge $\delta^{-2}\phi_1(a)\widetilde{\phi}_1(b)$}
\psfrag{del}{\huge $\delta$}
\resizebox{11.0cm}{!}{\includegraphics{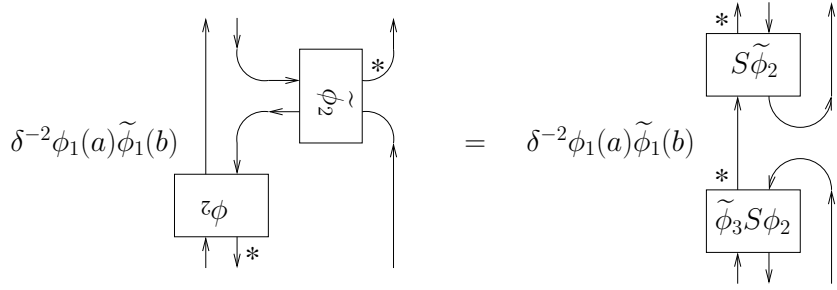}}
\end{center}
\caption{Simplifying the tangle on the left of Figure \ref{fig:needs}}
\label{fig:needs2}
\end{figure}

We now simplify:
\begin{eqnarray*}
\delta^{-2}\phi_1(a)\widetilde{\phi}_1(b)S\widetilde{\phi}_2 \otimes 
\widetilde{\phi}_3S\phi_2 &=& 
\delta^{-2}(\phi_1\widetilde{\phi}_3)(a)\widetilde{\phi}_1(b)S\widetilde{\phi}_2 
\otimes S\phi_2\\ 
&=& 
\delta^{-2}(\phi_1\widetilde{\phi}_1)(a)\widetilde{\phi}_2(b)S\widetilde{\phi}_3 
\otimes S\phi_2\\
&=& 
\delta^{-2}\phi_1(a_1)\widetilde{\phi}_1(a_2)\widetilde{\phi}_2(b)S\widetilde{\phi}_3 
\otimes S\phi_2\\
&=& 
\delta^{-2}\phi_1(a_1)\widetilde{\phi}_1(a_2b)S\widetilde{\phi}_2 \otimes 
S\phi_2\\
&=& \delta^{-1}\widetilde{\phi}_1(a_2b)S\widetilde{\phi}_2 \otimes 
\delta^{-1}\phi_1(a_1)S\phi_2\\
&=& SF(a_2b) \otimes SF(a_1)
\end{eqnarray*}
where the first equality is a consequence of $\phi_1 \otimes \psi S\phi_2 = 
\phi_1 \psi \otimes S\phi_2$ and the second is a consequence of the traciality 
of $\phi$.

Comparing the initial and terminal expressions in the above chain of 
equalities with the labelled tangles on the right in Figures \ref{fig:needs2} 
and \ref{fig:needs} completes the proof.
\end{proof}

We conclude with the following corollary.

\begin{Corollary}
Suppose that $N$ is a planar network with $g$ boxes all of which are $2$-boxes.
Then,
\[
Z^{P(H,\delta)}_N = 
Z^{P(H^*,\delta)}_{N^-} \circ F^{\otimes g},\]
where both sides are regarded as $\k$-valued functions on $H^{\otimes g}$.
\end{Corollary}

\begin{proof}
Note first that since $P(H,\delta)_{0_\pm}$ and $P(H^*,\delta)_{0_\pm}$ are
identified canonically with $\k$, the planar algebra isomorphism 
of Theorem \ref{theo:pliso}, which maps $^-P(H,\delta)_{0_\pm} =
P(H,\delta)_{0_\mp}$ to $P(H^*,\delta)_{0_\pm}$, is identified
with the identity map of $\k$. 

According to Theorem  \ref{theo:pliso}, if $N$ is a planar network with $g$ 
boxes all of which are $2$-boxes, then, with the identifications above, 
$Z^{^-P(H,\delta)}_N
= Z^{P(H^*,\delta)}_N \circ (SF)^{\otimes g}$.
But by definition, the former is $Z^{P(H,\delta)}_{N^-}$ while the latter is 
nothing but $Z^{P(H^*,\delta)}_{N^{--}} \circ (F)^{\otimes g}$, since (i) 
$N^{--}
= N \circ_{(B_1,\cdots ,B_g)}(R_2,\cdots ,R_2)$ where $R_2$ is the 
$2$-rotation tangle and (ii)
$Z^{P(H^*,\delta)}_{R_2} = S$.
Now, replacing $N^-$ by $N$ yields the desired conclusion.
\end{proof}

\section{appendix: tilings and tangles}

This brief appendix will be devoted to a statement of a result in 
combinatorial topology and a sketch of its application in proving
a generalisation of Lemma \ref{lemma:surj}. We omit all proofs.

Consider a convex $2k$-gon in the plane with its vertices numbered from
$1$ to $2k$ in clockwise order. By a tiling (by quadrilaterals) of the 
$2k$-gon we will
mean a collection of its diagonals that are required to be non-intersecting
and divide the polygon into quadrilaterals.

By a hexagon move on a tiling of a $2k$-gon, we will mean the following.
Take two of its quadrilaterals that share an edge and consider the hexagon
formed by their remaining edges. The common edge is a principal diagonal
of this hexagon. Replace this principal diagonal with one of the other two
principal diagonals of the hexagon to get a new tiling of the $2k$-gon.  

%The following proposition holds.
\begin{Proposition}\label{prop:tiling}
Any two tilings of a $2k$-gon are related
by a sequence of hexagon moves.\qed
\end{Proposition}

\begin{figure}[!h]
\begin{center}
%\psfrag{1}{\huge $1$}
%\psfrag{2}{\huge $2$}
%\psfrag{3}{\huge $3$}
%\psfrag{4}{\huge $4$}
%\psfrag{5}{\huge $5$}
%\psfrag{6}{\huge $6$}
\resizebox{11.0cm}{!}{\includegraphics{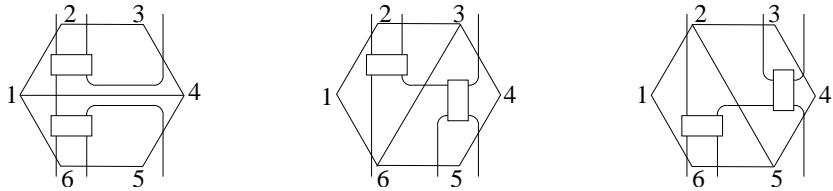}}
\end{center}
\caption{Some tangles and their associated tilings}
\label{fig:hexmoves}
\end{figure}

The point of this digression into tilings and hexagon moves is roughly that
tilings of a $2k$-gon correspond to tangles in ${\mathcal T(k)}$ (modulo
the equivalence relation that forgets the internal $*$'s)
while hexagon moves correspond to applying the 
exchange relations (E) and (A).
Figure \ref{fig:hexmoves} shows some elements of 
${\mathcal T}(3)$ 
and
their corresponding tilings of the hexagon.
Proposition \ref{prop:tiling} is the main step in the following
generalisation of Lemma \ref{lemma:surj}.

\begin{Lemma} For each $k \geq 3$, and each $X \in {\mathcal T}(k)$, the
map $Z^P_{X}
 : P_2^{\otimes (k-1)} \rightarrow P_k$ 
is surjective.\qed
\end{Lemma}


\begin{thebibliography}{amsalpha}

\bibitem[DasKdy]{DasKdy} Paramita Das and Vijay Kodiyalam, {\em Planar algebras
and the Ocneanu-Szymanski theorem}, Proc. AMS, {\bf 133}, (2005) 2751-2759.


\bibitem[Jns]{Jns} V. F. R. Jones, {\em Planar algebras I}, New
 Zealand J. of Math., to appear. e-print arXiv: math.QA/9909027

\bibitem[KdyLndSnd]{KdyLndSnd} Vijay Kodiyalam, Zeph Landau and V. S. Sunder,
{\em The planar algebra associated to a Kac algebra}, Proc. Indian Acad. of 
Sciences, {\bf 113}, (2003) 15-51.

\bibitem[KdySnd]{KdySnd} Vijay Kodiyalam and V. S. 
Sunder, {\em On Jones' planar algebras}, J. Knot theory and its ramifications, 
{\bf 13}, (2004) 219-247.
 
\bibitem[Lnd]{Lnd} Zeph Landau, {\em Exchange Relation Planar
    Algebras}, Proceedings of the Conference on Geometric and
    Combinatorial Group Theory, Part II (Haifa, 2000).  Geom. Dedicata
    {\bf 95}  (2002), 183--214.  

\bibitem[LndSnd]{LndSnd} Zeph Landau and V.S. Sunder, {\em Planar
    depth and planar subalgebras},  J. of Functional Analysis, {\bf
    195}, (2002), 71-88. 

\bibitem[LrsRdf]{LrsRdf} R.G. Larson and D. E. Radford, {\em Finite dimensional
cosemisimple Hopf algebras in characteristic $0$ are semisimple}, J. of 
Algebra,
{\bf 117}, 1988, 267--289.

\bibitem[Szy]{Szy} Wojciech Szymanski, {\em Finite index subfactors and
      Hopf algebra crossed products.} Proc. Amer. Math. Soc. {\bf 120}
      (1994), 519­528. 

\bibitem[TngGlk]{TngGlk} Pavel Etingof and Shlomo Gelaki,  {\em On
 finite-dimensional  semisimple and cosemisimple Hopf algebras in
 positive characteristic.} Internat. Math. Res. Notices  1998,
 {\bf 16}, 851--864.  


\end{thebibliography}
\end{document}